\documentclass{amsart}
\usepackage{amssymb}


\newcommand{\dbar}{\ensuremath{\overline\partial}}
\newcommand{\dbarstar}{\ensuremath{\overline\partial^*}}
\newcommand{\C}{\ensuremath{\mathbb{C}}}
\newcommand{\R}{\ensuremath{\mathbb{R}}}
\makeatletter     
\newcommand{\sumprime}{\if@display\sideset{}{'}\sum%
            \else\sum'\fi}    
\makeatother

\newtheorem{theorem}{Theorem}
\newtheorem{proposition}[theorem]{Proposition}
\newtheorem{lemma}[theorem]{Lemma}
\newtheorem{corollary}[theorem]{Corollary}

\DeclareMathOperator{\dom}{dom}

\begin{document}

\title[Global regularity of the $\dbar$-Neumann problem]{Global
  regularity of the $\dbar$-Neumann problem: a survey of the
  $L^2$-Sobolev theory}

\author{Harold P. Boas}
\address{Department of Mathematics\\
         Texas A\&M University\\
         College Station, TX 77843-3368}
\email{boas@math.tamu.edu}
\author{Emil J. Straube}
\email{straube@math.tamu.edu}
\thanks{Both authors partially supported by NSF grant number
  DMS-9500916 and at the Mathematical Sciences Research Institute 
by NSF grant number DMS-9022140.}

\subjclass{32F20, 32F15, 32H10, 35N15}

\maketitle

\tableofcontents

\section{Introduction}
The $\dbar$-Neumann problem is a natural example of a
boundary-value problem with an elliptic operator but with
non-coercive boundary conditions. It is also a prototype (in the
case of finite-type domains) of a subelliptic boundary-value
problem, in much the same way that the Dirichlet problem is the
archetypal elliptic boundary-value problem.  In this survey, we
discuss global regularity of the $\dbar$-Neumann problem in the
$L^2$-Sobolev spaces $W^s(\Omega)$ for all non-negative~$s$ and
also in the space $C^\infty(\overline\Omega)$.  For estimates in
other function spaces, such as H\"older spaces and $L^p$-Sobolev
spaces, see \cite{BealsGreinerStanton1987, BerndtssonB1994,
  ChangNagelStein1992, ChoS1995, ChristM1991, FeffermanC1995,
  FeffermanKohn1988, FeffermanKohnMachedon1990, GreinerStein1977,
  KerzmanN1971, KrantzS1979, LiebI1993, McNealJ1991,
  McNealStein1994, NagelRosaySteinWainger1989, StraubeE1995}; for
questions of real analytic regularity, see, for example,
\cite{ChenS1988, ChristM1996b, DerridjTartakoff1976,
  KomatsuG1976, TartakoffD1978, TartakoffD1980, TolliF1996,
  TrevesF1978} and section~10 of Christ's article
\cite{ChristM1997} in these proceedings.

We also discuss the closely related question of global regularity
of the Bergman projection operator. This question is intimately
connected with the boundary regularity of holomorphic mappings
(see, for example, \cite{BedfordE1984, BellS1981, BellS1984,
  BellS1990, BellLigocka1980,
  BellCatlin1982, DiederichFornaess1982, ForstnericF1993}). 

For an overview of techniques of partial differential equations
in complex analysis, see \cite{FollandKohn1972, HormanderL1965,
  HormanderL1990, HormanderL1994, KohnJ1977, KrantzS1992b}.

\section{The $L^2$ existence theory}\label{L2 theory}
Throughout the paper, $\Omega$ denotes a bounded domain
in~$\C^n$, where $n>1$.  We say that $\Omega$~has class~$C^k$
boundary if $\Omega=\{z: \rho(z)<0\}$, where $\rho$~is a
$k$~times continuously differentiable real-valued function in a
neighborhood of the closure~$\overline\Omega$ whose gradient is
normalized to length~$1$ on the boundary~$b\Omega$.  We denote
the standard $L^2$-Sobolev space of order~$s$ by $W^s(\Omega)$
(see, for example, \cite{AdamsR1975, LionsMagenes1972,
  TrevesF1975}).  The space of $(0,q)$ forms with coefficients in
$W^s(\Omega)$ is written $W^s_{(0,q)}(\Omega)$, the norm being
defined by
\begin{equation}
\bigl\|\sumprime_J a_J\,d\bar z_J \bigr\|_s^2 =
\sumprime_J \|a_J\|_s^2,
\end{equation}
where 
$d\bar z_J$ means $d\bar z_{j_1}\wedge d\bar z_{j_2}\wedge \dots
\wedge d\bar z_{j_q}$, and 
the prime indicates that the sum is taken over strictly
increasing $q$-tuples~$J$. 
We will consider the coefficients~$a_J$, originally
defined only for increasing multi-indices~$J$, to be defined for
other~$J$ so as to be antisymmetric functions of the indices.
For economy of notation, we
restrict attention to $(0,q)$ forms; modifications for $(p,q)$
forms are simple (because the $\dbar$~operator
does not see the $dz$ differentials). 

The $\dbar$ operator acts as usual on a $(0,q)$ form 
via
\begin{equation}
\dbar \biggl( \sumprime_J a_J \,d\bar z_J \biggr)
= \sum_{j=1}^n\sumprime_{J} \frac{\partial
a_J}{\partial \bar z_j} d\bar z_{jJ}.
\end{equation}
The domain of $\dbar: L^2_{(0,q)}(\Omega) \to
L^2_{(0,q+1)}(\Omega)$ consists of those forms~$u$ for which
$\dbar u$, defined in the sense of distributions, belongs to
$L^2_{(0,q+1)}(\Omega)$.  It is routine to check that $\dbar$~is
a closed, densely defined operator from $L^2_{(0,q)}(\Omega)$ to
$L^2_{(0,q+1)}(\Omega)$. Consequently, the Hilbert-space
adjoint~$\dbarstar$ also exists and defines a closed, densely
defined operator from $L^2_{(0,q+1)}(\Omega)$ to
$L^2_{(0,q)}(\Omega)$.

Suppose $u=\sumprime_J u_J\,d\bar z_J$ is continuously
differentiable on the closure~$\overline\Omega$, and $\psi$~is a
smooth test form. If the boundary~$b\Omega$ is sufficiently
smooth, then pairing $u$ with~$\dbar\psi$ and integrating by
parts gives
\begin{equation}
  (u,\dbar\psi)= \biggl( -\sum_{k=1}^n \sumprime_{K}
  \frac{\partial u_{kK}}{\partial z_k}\,d\bar z_K,
  \psi \biggr) + \sumprime_K \int_{b\Omega} \overline\psi_K
  \sum_{k=1}^n  u_{kK}
  \frac{\partial\rho}{\partial z_k}\,d\sigma.
\end{equation}
The same calculation with a compactly supported~$\psi$ shows
(without any boundary smoothness hypothesis) that if $u$~is a
square-integrable form in the domain of~$\dbarstar$, then
$\dbarstar u=\vartheta u$, where the \emph{formal
  adjoint~$\vartheta$} is given by the equation
\begin{equation}
  \vartheta u= -\sum_{k=1}^n \sumprime_{K}
 \frac{\partial u_{kK}}{\partial z_k}\,d\bar z_K.
\end{equation}
It follows that a continuously differentiable form~$u$ is in
the domain of~$\dbarstar$ if and only if
\begin{equation}
  \label{boundary condition}
  \sum_{k=1}^n
  u_{kK}\frac{\partial\rho}{\partial z_k}\biggr|_{b\Omega}=0 \text{
    for every~$K$}.
\end{equation}

The method of Friedrichs mollifiers shows that forms which are
continuously differentiable on the closure~$\overline\Omega$ are
dense in the intersection of the domains of $\dbar$
and~$\dbarstar$ with respect to the graph norm $(\|u\|^2 +\|\dbar
u\|^2 + \|\dbarstar u\|^2)^{1/2}$ when the boundary~$b\Omega$ is
sufficiently smooth (see, for instance, \cite[\S1.2 and
Prop.~2.1.1]{HormanderL1965}).  Also, forms that are continuously
differentiable on the closure are dense in the domain of~$\dbar$
with respect to the graph norm $(\|u\|^2+\|\dbar u\|^2)^{1/2}$.

The fundamental $L^2$ existence theorem for the $\dbar$-Neumann
problem is due to H\"ormander \cite{HormanderL1965}.
One version of the result is the following.

\begin{theorem}
  \label{L2 existence theorem}
  Let $\Omega$~be a bounded pseudoconvex domain in~$\C^n$, where
  $n\ge2$. Let $D$~denote the diameter of~$\Omega$, and suppose
  $1\le q\le n$.
  \begin{enumerate}
  \item The complex Laplacian $\square=\dbar \dbarstar +
    \dbarstar \dbar$ is an unbounded, self-adjoint, surjective
    operator from $L^2_{(0,q)}(\Omega)$ to itself having a
    bounded inverse~$N_q$ (the $\dbar$-Neumann operator).
  \item For all $u$ in $L^2_{(0,q)}(\Omega)$, we have the
    estimates
    \begin{equation}
    \begin{aligned}
      \|N_q u\|&\le \left(\frac{D^2e}{q}\right)\|u\| \\
      \|\dbarstar N_q u\|&\le \left(\frac{D^2e}{q}\right)^{1/2} \|u\| \\
      \|\dbar N_q u\|&\le \left(\frac{D^2e}{q}\right)^{1/2} \|u\|.
    \end{aligned}\label{L2 N estimate}
    \end{equation}
  \item If $f$~is a $\dbar$-closed $(0,q)$~form, then the
    canonical solution of the equation $\dbar u=f$ (the solution
    orthogonal to the kernel of~$\dbar$) is given by $u=\dbarstar
    N_q f$; if $f$~is a $\dbarstar$-closed $(0,q)$~form, then the
    canonical solution of the equation $\dbarstar u=f$ (the
    solution orthogonal to the kernel of~$\dbarstar$) is given by
    $u=\dbar N_q f$.
  \end{enumerate}
\end{theorem}

The Hilbert space method for proving Theorem~\ref{L2 existence
  theorem} is based on estimating the norm of a form~$u$ in terms
of the norms of $\dbar u$ and $\dbarstar u$.  H\"ormander
discovered that it is advantageous to introduce weighted spaces
$L^2(\Omega,e^{-\varphi})$, even for studying the unweighted
problem.  We denote the norm in the weighted space by
\(\|u\|_\varphi=\|u e^{-\varphi/2}\|\) and the adjoint of~$\dbar$
with respect to the weighted inner product by
\(\dbarstar_\varphi(\,\cdot\,)=e^\varphi\,
\dbarstar(\,\cdot\,e^{-\varphi} )\). More generally, one can
choose different exponential weights in $L^2_{(0,q-1)}$,
$L^2_{(0,q)}$, and $L^2_{(0,q+1)}$; see \cite{HormanderL1990} for
this method and applications.

The following identity is the basic starting point. The proof
involves integrating by parts and manipulating the boundary
integrals with the aid of the boundary condition~\eqref{boundary
  condition} for membership in the domain of~$\dbarstar$.  The
idea of introducing a second auxiliary function~$a$ originated
with Ohsawa and Takegoshi \cite{OhsawaTakegoshi1987, OhsawaT1988}
in their work on extending square-integrable holomorphic
functions from submanifolds.  The formulation given below comes
from the recent work of Siu \cite{SiuY1996} and McNeal
\cite{McNealJ1996}. In these papers (see also
\cite{DiederichHerbort1990}), the freedom to manipulate both the
weight factor~\(\varphi\) and the twisting factor~\(a\) is
essential.

\begin{proposition}
\label{prop:identity}
  Let $\Omega$~be a bounded domain in~$\C^n$ with class~$C^2$
  boundary; let $u$~be a $(0,q)$ form (where $1\le q\le n$) that
  is in the domain of~$\dbarstar$ and that is continuously
  differentiable on the closure~$\overline\Omega$; and let $a$
  and~$\varphi$ be real functions that are twice continuously
  differentiable on~$\overline\Omega$, with $a\ge0$. Then
\begin{multline}
  \|\sqrt a\,\dbar u\|^2_\varphi + \|\sqrt a\,\dbarstar_\varphi
  u\|^2_\varphi = \sumprime_K \sum_{j,k=1}^n \int_{b\Omega} a
  \frac{\partial^2 \rho}{\partial z_j \partial\bar z_k} u_{jK}
  \bar u_{kK} e^{-\varphi}\, d\sigma \\ +\sumprime_J \sum_{j=1}^n
  \int_\Omega a \left| \frac{\partial u_J}{\partial \bar
      z_j}\right|^2 e^{-\varphi}\, dV + 2\Re \biggl( \sumprime_K
  \sum_{j=1}^n u_{jK} \frac{\partial a}{\partial z_j} \, d\bar
  z_K, \dbarstar_\varphi u\biggr)_\varphi \\ +\sumprime_K
  \sum_{j,k=1}^n \int_{\Omega} \left( a
    \frac{\partial^2\varphi}{\partial z_j\partial \bar z_k}
    -\frac{\partial^2 a}{\partial z_j \partial \bar z_k} \right)
  u_{jK} \bar u_{kK} e^{-\varphi}\, dV.
\label{eqn:identity}
\end{multline}
\end{proposition}

For \(a\equiv1\) see \cite{HormanderL1965}; the case \(a\equiv1\)
and \(\varphi\equiv0\) is the classical Kohn-Morrey inequality
\cite{KohnJ1963, KohnJ1964, MorreyC1958} (see also
\cite{AshM1964}).  The usual proof of the $L^2$~existence theorem
is based on a variant of~\eqref{eqn:identity} with $a\equiv1$ and
with different exponential weights~$\varphi$ in the different
\(L^2_{(0,q)}\) spaces; see \cite{CatlinD1984b} for an elegant
implementation of this approach. Here we will give an argument
that has not appeared explicitly in the literature: we take
$\varphi\equiv0$ and make a good choice of~$a$.

Suppose that $\Omega$~is a pseudoconvex domain: this means
that the complex Hessian of the defining function~$\rho$ is a
non-negative form on the vectors in the complex tangent space.
Consequently, the boundary integral in~\eqref{eqn:identity} is
non-negative.  In particular, taking~$a$ to be identically
equal to~$1$ gives that
\begin{equation}
\|\dbar u\|^2+ \|\dbarstar
u\|^2 \ge \sumprime_J\sum_{j=1}^n \|\partial u_J/\partial \bar
z_j\|^2, \label{free}
\end{equation}
so the bar derivatives of~$u$ are always under control.

If we replace~$a$ by $1-e^b$, where $b$~is an arbitrary twice
continuously differentiable non-positive function, then after
applying the Cauchy-Schwarz inequality to the term
in~\eqref{eqn:identity} involving first derivatives of~$a$, we
find
\begin{equation}
  \|\sqrt a\,\dbar u\|^2 + \|\sqrt a\,\dbarstar u\|^2 \ge
  \sumprime_K \sum_{j,k=1}^n \int_\Omega e^b \frac{\partial^2
    b}{\partial z_j \partial \bar z_k} u_{jK} \bar u_{kK} \, dV
  -\|e^{b/2}\dbarstar u\|^2.
\end{equation}
Since $a+e^b=1$ and $a\le 1$, it follows that
\begin{equation}
\label{basic inequality}
\|\dbar u\|^2 + \|\dbarstar u\|^2 \ge \sumprime_K \sum_{j,k=1}^n
\int_\Omega e^b \frac{\partial^2 b}{\partial z_j \partial \bar
  z_k} u_{jK} \bar u_{kK} \, dV
\end{equation}
for every twice continuously differentiable non-positive
function~$b$. Notice that this inequality becomes a strong one if
there happens to exist a bounded plurisubharmonic function~$b$
whose complex Hessian has large eigenvalues. (This theme will
recur later on: see the discussion after
Theorem~\ref{thm:subelliptic} and the discussion of property~(P)
in section~\ref{sec:compactness}.)

In particular, let $p$~be a point of~$\Omega$, and set $b(z)=-1 +
|z-p|^2/D^2$, where $D$~is the diameter of the bounded
domain~$\Omega$.  The preceding inequality then implies the
fundamental estimate
\begin{equation}
\|u\|^2 \le \frac{D^2 e}{q}(\|\dbar u\|^2 + \|\dbarstar u\|^2).
\label{L2 estimate}
\end{equation}
Although this estimate was derived under the assumption that
$u$~is continuously differentiable on the
closure~$\overline\Omega$, it holds by density for all
square-integrable forms~$u$ that are in the intersections of the
domains of $\dbar$ and~$\dbarstar$.  We also assumed that the
boundary of~$\Omega$ is smooth enough to permit integration by
parts.  
Estimate~\eqref{L2 estimate} is equivalent to every form in
\(L^2_{(0,q)}(\Omega)\) admitting a representation as \(\dbar
v+\dbarstar w\) with \(\|v\|^2+\|w\|^2\le (D^2 e/q)
\|u\|^2\). The latter property carries over to arbitrary bounded
pseudoconvex domains by exhausting a nonsmooth~$\Omega$ by smooth
ones, and therefore so does inequality~\eqref{L2 estimate}.

Once estimate~\eqref{L2 estimate} is in hand, the proof of
Theorem~\ref{L2 existence theorem} follows from standard Hilbert
space arguments; see, for example,
\cite[pp.~164--165]{CatlinD1983} or \cite[\S2]{ShawM1992}. The
latter paper also shows the existence of the $\dbar$-Neumann
operator~$N_0$ on $(\ker\dbar)^\perp$.

\section{Regularity on general pseudoconvex domains}\label{general}
A basic question is whether one can improve Theorem~\ref{L2
  existence theorem} to get regularity estimates in Sobolev
norms: $\|Nu\|_s\le C\|u\|_s$, $\|\dbarstar Nu\|_s\le C\|u\|_s$,
$\| \dbar N u\|_s\le C\|u\|_s$. If such estimates were to hold
for all positive~$s$, then Sobolev's lemma would imply that the
$\dbar$-Neumann operator~$N$ (together with $\dbarstar N$ and
$\dbar N$) is continuous in the space $C^\infty(\overline\Omega)$
of functions smooth up to the boundary.

At first sight, it appears that one ought to able to generalize
the fundamental $L^2$ estimate~\eqref{L2 estimate} directly to an
estimate of the form $\|u\|_s\le C(\|\dbar u\|_s+ \|\dbarstar
u\|_s)$, simply by replacing~$u$ by a derivative of~$u$. This
naive expectation is erroneous: the difficulty is that not every
derivative of a form~$u$ in the domain of~$\dbarstar$ is again in
the domain of~$\dbarstar$.  The usual attempt to overcome this
difficulty is to cover the boundary of~$\Omega$ with special
boundary charts \cite[p.~33]{FollandKohn1972} in each of which
one can take a frame of tangential vector fields that do preserve
the domain of~$\dbarstar$.  Since such vector fields have
variable coefficients, they do not commute with either $\dbar$
or~$\dbarstar$, and so one needs to handle error terms that arise
from the commutators.

In subsequent sections, we will discuss various hypotheses on the
domain~$\Omega$ that yield regularity estimates in Sobolev norms.
In this section, we discuss firstly some completely general
results on smoothly bounded pseudoconvex domains and secondly
some counterexamples.

It is an observation of J.~J.~Kohn and his school that the
$\dbar$-Neumann problem is always regular in $W^\epsilon(\Omega)$
for a sufficiently small positive~$\epsilon$.

\begin{proposition}
\label{thm:epsilon}
  Let $\Omega$ be a bounded pseudoconvex domain in~$\C^n$ with
  class~$C^\infty$ boundary. 
  There exist positive $\epsilon$ and~$C$
  (both depending on~$\Omega$) such that $\|Nu\|_\epsilon\le
  C\|u\|_\epsilon$, $\|\dbarstar Nu\|_\epsilon \le
  C\|u\|_\epsilon$, and $\|\dbar Nu\|_\epsilon\le
  C\|u\|_\epsilon$ for every $(0,q)$ form~$u$ (where $1\le q\le
  n$).
\end{proposition}

A proof seems never to have appeared in print, but the idea is
very simple.  Since the commutator of a differential operator of
order~$\epsilon$ with $\dbar$ or~$\dbarstar$ is again an operator
of order~\(\epsilon\), but with a coefficient bounded by a
constant times~$\epsilon$, error terms can be absorbed into the
main term when $\epsilon$~is sufficiently small.  

\begin{theorem}
  \label{global regularity}
  Let $\Omega$ be a bounded pseudoconvex domain in~$\C^n$ with
  class~$C^\infty$ boundary. Fix a positive~$s$. There exists
  a~$T$ (depending on $s$ and~$\Omega$) such that for every~$t$
  larger than~$T$, the weighted $\dbar$-Neumann problem for the
  space $L^2_{(0,q)}(\Omega, e^{-t|z|^2}\,dV(z))$ is regular in
  $W^s(\Omega)$. In other words, $N_t$, $\dbarstar_t N_t$, and
  $\dbar N_t$ are continuous in $W^s(\Omega)$.

  Moreover, if $f$~is a $\dbar$-closed $(0,q)$ form with
  coefficients in $C^\infty(\overline\Omega)$, then there exists
  a form~$u$ with coefficients in $C^\infty(\overline\Omega)$ such
  that $\dbar u=f$.

\end{theorem}

This fundamental result on continuity of the weighted operators
is due to Kohn \cite{KohnJ1973}. It says that one can always have
regularity for the $\dbar$-Neumann problem up to a certain number
of derivatives if one is willing to change the measure with
respect to which the problem is defined.  The idea of the proof
is to apply Proposition~\ref{prop:identity} with \(a\equiv1\) and
\(\varphi(z)=t|z|^2\) to obtain $\|e^{-t|z|^2/2}u\|^2 \le
Ct^{-1}(\|e^{-t|z|^2/2}\,\dbar u\|^2
+\|e^{-t|z|^2/2}\,\dbarstar_t u\|^2)$. When $t$~is sufficiently
large, the factor $t^{-1}$ makes it possible to absorb error
terms coming from commutators (see the sketch of the proof of
Theorem~\ref{compactness} below for the ideas of the technique).
The resulting \emph{a priori} estimates are valid under the
assumption that the left-hand sides of the inequalities are known
to be finite; Kohn completed the proof by applying the method of
elliptic regularization \cite{KohnNirenberg1965} (see also the
remarks after Theorem~\ref{strictly} below).

Via a Mittag-Leffler argument (\cite[p.~230]{KohnJ1977}, argument
attributed to H\"ormander), one can deduce solvability of the
equation $\dbar u=f$ in the space $C^\infty(\overline\Omega)$
(but the solution will not be the canonical solution orthogonal
to the kernel of~$\dbar$).  With some extra care, the solution
operator can be made linear, and also continuous from
$W^{s+\epsilon}_{(0,q+1)}(\Omega)\cap\ker\dbar$ to
$W^s_{(0,q)}(\Omega)$ for every positive $s$ and~$\epsilon$
\cite{SibonyN}. It is unknown whether or not there exists a
linear solution operator for~$\dbar$ that breaks even at every
level in the Sobolev scale.  Solvability with Sobolev estimates
(with a loss of three derivatives) has recently been obtained for
domains with only $C^4$~boundary by S.~L.~Yie \cite{YieS}.

Given any solution of the equation $\dbar u=f$, one obtains
the canonical solution by subtracting from~$u$ its projection
onto the kernel of~$\dbar$.  In view of Kohn's result above,
it is natural to study the regularity properties of the
projection mapping.  We denote the orthogonal projection from
$L^2_{(0,q)}(\Omega)$ onto $\ker \dbar$ by~$P_q$; when $q=0$,
this operator is the \emph{Bergman projection}. A direct
relation between the Bergman projection and the
$\dbar$-Neumann operator is given by Kohn's formula
$P_q=\mathrm{Id}- \dbarstar N_{q+1}\dbar$ for $0\le q\le
n$. It is evident that if the $\dbar$-Neumann operator
$N_{q+1}$ is continuous in $C^\infty(\overline\Omega)$, then
so is~$P_q$. The 
exact relationship between regularity properties of the
$\dbar$-Neumann operators and the Bergman projections was
determined in
\cite{BoasStraube1990}.
\begin{theorem}
  \label{equivalence}
  Let $\Omega$ be a bounded pseudoconvex domain in~$\C^n$ with
  class~$C^\infty$ boundary. Fix an integer~$q$ such that $1\le
  q\le n$. Then the $\dbar$-Neumann operator~$N_q$ is continuous
  on $C^\infty_{(0,q)}(\overline\Omega)$ if and only if the
  projection operators $P_{q-1}$, $P_q$, and $P_{q+1}$ are
  continuous on the corresponding $C^\infty(\overline\Omega)$
  spaces. The analogous statement holds with the Sobolev space
  $W^s(\Omega)$ in place of~$C^\infty(\overline\Omega)$.
\end{theorem}

In view of the implications for boundary regularity of
biholomorphic and proper holomorphic mappings \cite{BedfordE1984,
  BellS1981, BellS1984, BellS1990, BellLigocka1980,
  BellCatlin1982, DiederichFornaess1982, ForstnericF1993},
regularity in \(C^\infty(\overline\Omega)\) is a key issue.

For some years there was uncertainty over whether the Bergman
projection operator~\(P_0\) of every bounded domain in~$\C^n$
with $C^\infty$~smooth boundary might be regular in the space
$C^\infty(\overline\Omega)$. Barrett \cite{BarrettD1984} found
the first counterexample, motivated by the so-called ``worm
domains'' of Diederich and Forn{\ae}ss
\cite{DiederichFornaess1977a}. In his example, for every $p>2$
there is a smooth, compactly supported function whose Bergman
projection is not in $L^p(\Omega)$.
In~\cite{BarrettFornaess1986}, Barrett and Forn{\ae}ss
constructed a counterexample even more closely related to the
worm domains.  Although the worm domains are smoothly bounded
pseudoconvex domains in~$\C^2$, these counterexamples are not
pseudoconvex.  Subsequently, Kiselman \cite{KiselmanC1991} showed
that pseudoconvex, but nonsmooth, truncated versions of the worm
domains have irregular Bergman projections.

Later Barrett \cite{BarrettD1992} (see \cite{BarrettD1997} for a
generalization) used a scaling argument together with
computations on piecewise Levi-flat model domains to show that
the Bergman projection of a worm domain must fail to preserve the
space $W^s(\Omega)$ when $s$~is sufficiently large.  In view of
Theorem~\ref{equivalence}, the $\dbar$-Neumann operator~$N_1$
also fails to preserve $W^s_{(0,1)}(\Omega)$.  This result left
open the possibility of regularity in
$C^\infty(\overline\Omega)$.  Finally the question was resolved
by Christ \cite{ChristM1996a}, as follows.
\begin{theorem}
  \label{worm}
  For every worm domain, the Bergman projection operator~$P_0$
  and the $\dbar$-Neumann operator~$N_1$ fail to be continuous on
  $C^\infty(\overline\Omega)$ and
  $C^\infty_{(0,1)}(\overline\Omega)$.
\end{theorem}
Christ's proof is delicate and indirect. Roughly speaking, he
shows that the $\dbar$-Neumann operator does satisfy for most
values of~$s$ an estimate of the form $\|N_1 u\|_s\le C\|u\|_s$
for all~$u$ for which $N_1 u$ is known a priori to lie in
$C^\infty_{(0,1)}(\overline\Omega)$. If $N_1$ were to preserve
$C^\infty_{(0,1)}(\overline\Omega)$, then density of
$C^\infty_{(0,1)}(\overline\Omega)$ in $W^s_{(0,1)}(\Omega)$
would imply continuity of~$N_1$ in $W^s_{(0,1)}(\Omega)$,
contradicting Barrett's result.

The obstruction to continuity in \(W^s(\Omega)\) for every~\(s\)
on the worm domains is a global one: namely, the nonvanishing of
a certain class in the first De~Rham cohomology of the annulus of
weakly pseudoconvex boundary points (this class measures the
twisting of the boundary at the annulus; for details, see
Theorem~\ref{cohomology}).  For smoothly bounded
domains~\(\Omega\), it is known that for each fixed~\(s\) there
is no local obstruction in the boundary to continuity in
\(W^s(\Omega)\) \cite{BarrettD1986, ChenS1991b}.

For all domains~\(\Omega\) where continuity in
\(C^\infty(\overline\Omega)\) is known, one can actually prove
continuity in \(W^s(\Omega)\) for all positive~\(s\).  This
intriguing phenomenon is not understood at present.  (The
corresponding phenomenon does not hold for partial differential
operators in general: see the discussion in section~3 of Christ's
article \cite{ChristM1997} in these proceedings.)

Although regularity of the $\dbar$-Neumann problem in
$C^\infty(\overline\Omega)$ is known in large classes of
pseudoconvex domains (see sections \ref{sec:finite
  type}--\ref{vector field}), the example of the worm domains
shows that regularity sometimes fails.  At present, necessary and
sufficient conditions for global regularity of the
$\dbar$-Neumann operator and of the Bergman projection are not
known.

\section{Domains of finite type}
\label{sec:finite type}
Historically, the first major development on the $\dbar$-Neumann
problem was its solution by Kohn \cite{KohnJ1963, KohnJ1964} for
strictly pseudoconvex domains. A strictly pseudoconvex domain can
be defined by a strictly plurisubharmonic function, so by taking
\(a\equiv1\) and \(\varphi\equiv0\) in~\eqref{eqn:identity} and
keeping the boundary term we find that $\|\dbar u\|^2+
\|\dbarstar u\|^2\ge C\|u\|^2_{L^2(b\Omega)}$. Roughly speaking,
this inequality says that we have gained half a derivative, since
the restriction map $W^{s+\frac12}(\Omega)\to W^s(b\Omega)$ is
continuous when $s>0$. This gain is half of what occurs for an
ordinary elliptic boundary-value problem, so we have a
``subelliptic estimate.''
\begin{theorem}
\label{strictly}
Let $\Omega$ be a bounded strictly pseudoconvex domain in~$\C^n$
with class $C^\infty$ boundary. If $1\le q\le n$, then for each
non-negative~$s$ there is a constant~$C$ such that the following
estimates hold for every $(0,q)$ form~$u$:
\begin{equation}
\begin{aligned}
\|u\|_{s+\frac12}&\le C(\|\dbar u\|_s+\|\dbarstar
u\|_s) \text{ if $u\in\dom\dbar\cap\dom\dbarstar$},\\
\|N_q u\|_{s+1}&\le C\|u\|_s,\\
\|\dbar N_q u\|_{s+\frac12}
+\|\dbarstar N_q u\|_{s+\frac12}&\le C\|u\|_s.
\end{aligned}
\end{equation}
\end{theorem}

The standard reference for the proof of this result is
\cite{FollandKohn1972} (where the theory is developed for almost
complex manifolds); see also \cite{KrantzS1992b}. The estimates
can be localized, as in Theorem~\ref{pseudolocal estimates}
below.  

A key technical point in the proof of Theorem~\ref{strictly} is
that after establishing the estimates under the assumption that
the left-hand side is a priori finite, one then has to convert
the a priori estimates into genuine estimates, in the sense that
the left-hand side is finite when the right-hand side is finite.
Kohn's original approach was considerably simplified in
\cite{KohnNirenberg1965} in a very general framework, via the
elegant device of ``elliptic regularization.'' The idea of the
method is to add to~$\square$ an elliptic operator
times~$\epsilon$ (thereby obtaining a standard elliptic problem),
to prove estimates independent of~$\epsilon$, and to let
$\epsilon$~go to zero.  (The analysis Christ used
\cite{ChristM1996a} to prove Theorem~\ref{worm} shows that indeed
a priori estimates cannot always be converted into genuine
estimates. For this phenomenon in the context of the Bergman
projection, see \cite{BoasStraube1992a}.)
Another interesting approach to the proof of
Theorem~\ref{strictly} was indicated by Morrey~\cite{MorreyC1963}.

A number of authors (see \cite{BealsGreinerStanton1987,
  GreinerStein1977} and their references) have refined the
results for strictly pseudoconvex domains in various ways, such
as estimates in other function spaces and anisotropic estimates.
In particular, \(N\)~gains two derivatives in complex tangential
directions; this gain results from the bar derivatives always
being under control (see~\eqref{free}).  Integral kernel methods
have also been developed successfully on strictly pseudoconvex
domains; see \cite{GrauertLieb1970, HenkinG1969, HenkinG1970,
  LiebI1993, LiebRange1987, RamirezE1970, RangeR1986, RangeR1987}
and their references.

The gain of one derivative for the $\dbar$-Neumann operator~$N_1$
in Theorem~\ref{strictly} is sharp, and the domain is necessarily
strictly pseudoconvex if this estimate holds. For discussion of
this point, see \cite{CatlinD1983},
\cite[\S{}III.2]{FollandKohn1972}, \cite[\S3.2]{HormanderL1965},
and~\cite[\S4]{KrantzS1979}.

More generally, one can ask when the $\dbar$-Neumann operator
gains some fractional derivative. One says that a subelliptic
estimate of order~$\epsilon$ holds for the $\dbar$-Neumann
problem on $(0,q)$ forms in a neighborhood~$U$ of a boundary
point~$z_0$ of a pseudoconvex domain in~$\C^n$ if there is a
constant~$C$ such that
\begin{equation}
\label{subelliptic estimate}
\|u\|_\epsilon^2\le C(\|\dbar u\|_0^2 +\|\dbarstar u\|_0^2)
\end{equation}
for every smooth $(0,q)$ form~$u$ that is supported in
$U\cap\overline\Omega$ and that is in the domain of~$\dbarstar$.
The systematic study of subelliptic estimates in
\cite{KohnNirenberg1965} provides the following ``pseudolocal
estimates.''
\begin{theorem}
  \label{pseudolocal estimates}
  Let $\Omega$ be a bounded pseudoconvex domain in~$\C^n$ with
  class~$C^\infty$ boundary. Suppose that a subelliptic
  estimate~\eqref{subelliptic estimate} holds in a
  neighborhood~$U$ of a boundary point~$z_0$. Let $\chi_1$ and
  $\chi_2$ be smooth cutoff functions supported in~$U$ with
  $\chi_2$ identically equal to~$1$ in a neighborhood of the
  support of~$\chi_1$. For every non-negative~$s$, there is a
  constant~$C$ such that the $\dbar$-Neumann operator~$N_q$ and
  the Bergman projection~$P_q$ satisfy the estimates
\begin{equation}
\begin{aligned}
  \|\chi_1 N_q u\|_{s+2\epsilon} &\le C(\|\chi_2 u\|_s +\|u\|_0),
  \qquad 1\le q\le n,\\ \|\chi_1 \dbarstar N_q u\|_{s+\epsilon}
  +\|\chi_1 \dbar N_q u\|_{s+\epsilon} &\le C(\|\chi_2
  u\|_s+\|u\|_0), \qquad 1\le q\le n,\\ \|\chi_1 P_q u\|_s &\le
  C(\|\chi_2 u\|_s +\|u\|_0),\qquad 0\le q\le n.
\end{aligned}
\end{equation}
\end{theorem}
Consequently, if a subelliptic estimate~\eqref{subelliptic
  estimate} holds in a neighborhood of every boundary point of a
smooth bounded pseudoconvex domain~$\Omega$ in~$\C^n$, then the
Bergman projection is continuous from $W^s_{(0,q)}(\Omega)$ to
itself, and the $\dbar$-Neumann operator is continuous from
$W^s_{(0,q)}(\Omega)$ to $W^{s+2\epsilon}_{(0,q)}(\Omega)$.

In a sequence of papers \cite{CatlinD1983, CatlinD1984a,
  CatlinD1987b, D'AngeloJ1979, D'AngeloJ1980, D'AngeloJ1982},
Kohn's students David Catlin and John D'Angelo resolved the
question of when subelliptic estimates hold in a neighborhood of
a boundary point of a smooth bounded pseudoconvex domain
in~$\C^n$.  The necessary and sufficient condition is that the
point have ``finite type'' in an appropriate sense. We briefly
sketch this work; for details, consult the above papers as well
as \cite{CatlinD1987a, D'AngeloJ1993, D'AngeloJ1995,
  DiederichLieb1981, GreinerP1974, KohnJ1979a, KohnJ1979b,
  KohnJ1981, KohnJ1984} and the survey by D'Angelo and Kohn
\cite{D'AngeloKohn1997} in these proceedings.

The simplest obstruction to a subelliptic estimate is the
presence of a germ of an analytic variety in the boundary of a
domain.  Indeed, examples show that local regularity of the
$\dbar$-Neumann problem fails when there are complex varieties in
the boundary; see \cite{CatlinD1981, DiederichPflug1981}.  If the
boundary is real-analytic near a point, then the absence of germs
of \(q\)-dimensional complex-analytic varieties in the boundary
near the point is necessary and sufficient for the existence of a
subelliptic estimate on \((0,q)\)-forms \cite{KohnJ1979b}. This
was first proved by combining a sufficient condition from Kohn's
theory of ideals of subelliptic multipliers \cite{KohnJ1979b}
with a theorem of Diederich and Forn{\ae}ss
\cite{DiederichFornaess1978} on analytic varieties.  Moreover,
Diederich and Forn{\ae}ss showed that a compact real-analytic
manifold contains no germs of complex-analytic varieties of
positive dimension, so subelliptic estimates hold for every
bounded pseudoconvex domain in~$\C^n$ with real-analytic
boundary.

The first positive results in the $C^\infty$~category were
established in dimension two. A boundary point of a domain
in~$\C^2$ is of finite type if the boundary has finite order of
contact with complex manifolds through the point; equivalently,
if some finite-order commutator of complex tangential vector
fields has a component that is transverse to the complex tangent
space to the boundary. If $m$~is an upper bound for the order of
contact of complex manifolds with the boundary, then a
subelliptic estimate \eqref{subelliptic estimate} holds with
$\epsilon=1/m$. For these results, see \cite{GreinerP1974,
  KohnJ1979b}; for the equivalence of the two notions of finite
type, see \cite{BloomGraham1977}. For pseudoconvex domains of
finite type in dimension two, sharp estimates for the
\dbar-Neumann problem are now known in many function spaces (see
\cite{ChangNagelStein1992, ChristM1991} and their references).

In higher dimensions, it is no longer the case that all
reasonable notions of finite type agree; for relations among
them, see \cite{D'AngeloJ1987a}.  D'Angelo's notion of finite
type has turned out to be the right one for characterizing
subelliptic estimates for the $\dbar$-Neumann problem.  His idea
to measure the order of contact of varieties with a real
hypersurface~$M$ in~$\C^n$ at a point~$z_0$ is to fix a defining
function~$\rho$ for~$M$ and to consider the order of vanishing at
the origin of $\rho\circ f$, where $f$~is a nonconstant
holomorphic mapping from a neighborhood of the origin in~$\C$
to~$\C^n$ with $f(0)=z_0$.  Since the variety that is the image
of~$f$ may be singular, it is necessary to normalize by dividing
by the order of vanishing at the origin of $f(\,\cdot\,)-z_0$.
The supremum over all~$f$ of this normalized order of contact of
germs of varieties with~$M$ is the D'Angelo $1$-type of~$z_0$.

\begin{theorem}
  The set of points of finite $1$-type of a smooth real
  hypersurface~$M$ in~$\C^n$ is an open subset of~$M$, and the
  $1$-type is a locally bounded function on~$M$.
\end{theorem}
This fundamental result of D'Angelo \cite{D'AngeloJ1982} is
remarkable, because the $1$-type may fail to be an upper
semi-continuous function (see \cite[p.~136]{D'AngeloJ1993} for a
simple example).  The theorem implies that if every point of a
bounded domain in~$\C^n$ is of finite $1$-type, then there is a
global upper bound on the $1$-type.

For higher-dimensional varieties, there is no canonical way that
serves all purposes to define the order of contact with a
hypersurface. Catlin \cite{CatlinD1987b} defined a quantity
$D_q(z_0)$ that measures the order of contact of $q$-dimensional
varieties in ``generic'' directions (and $D_1$~agrees with
D'Angelo's $1$-type).  Catlin's fundamental result is the
following.
\begin{theorem}
\label{thm:subelliptic}
  Let $\Omega$ be a bounded pseudoconvex domain in~$\C^n$ with
  class~$C^\infty$ boundary. A subelliptic estimate for the
  $\dbar$-Neumann problem on $(0,q)$ forms holds in a
  neighborhood of a boundary point~$z_0$ if and only if
  $D_q(z_0)$ is finite.  The $\epsilon$ in the subelliptic
  estimate \eqref{subelliptic estimate} satisfies $\epsilon \le
  1/D_q(z_0)$.
\end{theorem}

Catlin proved the necessity of finite order of contact, together
with the upper bound on~$\epsilon$, in \cite{CatlinD1983} (see
also \cite{CatlinD1981}), and the sufficiency in
\cite{CatlinD1987b}. Catlin's proof of sufficiency has two parts.
His theory of multitypes \cite{CatlinD1984a} implies the
existence of a stratification of the set of weakly pseudoconvex
boundary points. The stratification is used to construct families
of bounded plurisubharmonic functions whose complex Hessians in
neighborhoods of the boundary have eigenvalues that blow up like
inverse powers of the thickness of the neighborhoods.  Such powers
heuristically act like derivatives, and so it should be plausible
that the basic inequality~\eqref{basic inequality} leads to a
subelliptic estimate~\eqref{subelliptic estimate}.

It is unknown in general how to determine the optimal value
of~$\epsilon$ in a subelliptic estimate in terms of boundary
data. For convex domains of finite type in~$\C^n$, the
optimal~$\epsilon$ in a subelliptic estimate for $(0,1)$ forms is
the reciprocal of the D'Angelo $1$-type \cite{FornaessSibony1989,
  McNealJ1992}; this is shown by a direct construction of bounded
plurisubharmonic functions with suitable Hessians near the
boundary.  McNeal proved \cite{McNealJ1992} that for convex
domains, the D'Angelo $1$-type can be computed simply as the
maximal order of contact of the boundary with complex lines.
(There is an elementary geometric proof of McNeal's result
in \cite{BoasStraube1992b} and an analogue for Reinhardt
domains in \cite{FuKrantz}.) It is clear that in general, the
best~$\epsilon$ cannot equal the reciprocal of the type, simply
because the type is not necessarily upper semi-continuous. For
more about this subtle issue, see \cite{D'AngeloJ1993,
  D'AngeloJ1995, D'AngeloKohn1997, DiederichHerbort1993}.

\section{Compactness}
\label{sec:compactness}
A subelliptic estimate~\eqref{subelliptic estimate} implies, in
particular, that the $\dbar$-Neumann operator is compact as an
operator from $L^2_{(0,q)}(\Omega)$ to itself. This follows
because the embedding from $W^\epsilon_{(0,q)}(\Omega)$ into
$L^2_{(0,q)}(\Omega)$ is compact when $\Omega$~is bounded with
reasonable boundary, by the Rellich-Kondrashov theorem (see, for
example, \cite[Prop.~25.5]{TrevesF1975}). One might think of
compactness in the $\dbar$-Neumann problem as a limiting case of
subellipticity as $\epsilon\to0$.

The following lemma reformulates the compactness condition.

\begin{lemma}
  Let $\Omega$ be a bounded pseudoconvex domain in~$\C^n$, and
  suppose that $1\le q\le n$. The following statements are
  equivalent.
  \begin{enumerate}
  \item The $\dbar$-Neumann operator~$N_q$ is compact from
    $L^2_{(0,q)}(\Omega)$ to itself.\label{compactness:a}
  \item The embedding of the space $\dom\dbar\cap\dom\dbarstar$,
    provided with the graph norm $u\mapsto \|\dbar
    u\|_0+\|\dbarstar u\|_0$, into $L^2_{(0,q)}(\Omega)$ is
    compact.\label{compactness:b}
  \item For every positive~$\epsilon$ there exists a
    constant~$C_\epsilon$ such that
    \begin{equation}
      \|u\|_0^2\le\epsilon(\|\dbar u\|_0^2+\|\dbarstar u\|_0^2)
      +C_\epsilon\|u\|_{-1}^2 \label{compactness estimate}
    \end{equation}
  when $u\in\dom\dbar\cap\dom\dbarstar$.\label{compactness:c}
  \end{enumerate}
\end{lemma}

Statement~\ref{compactness:c} is called a \emph{compactness
  estimate} for the $\dbar$-Neumann problem. Its equivalence with
statement~\ref{compactness:b} is in
\cite[Lemma~1.1]{KohnNirenberg1965}. The equivalence of
statement~\ref{compactness:a} with statements \ref{compactness:b}
and~\ref{compactness:c} follows easily from the $L^2$~theory
discussed in section~\ref{L2 theory} and the compactness of the
embedding $L^2_{(0,q)}(\Omega)\to W^{-1}_{(0,q)}(\Omega)$.

In view of Theorem~\ref{pseudolocal estimates}, it is a
reasonable guess that compactness in the $\dbar$-Neumann problem
implies global regularity of the $\dbar$-Neumann operator in the
sense that $N_q$~maps $W^s_{(0,q)}(\Omega)$ into itself. Work of
Kohn and Nirenberg \cite{KohnNirenberg1965} shows that
this conjecture is correct. 

\begin{theorem}
  \label{compactness}
  Let $\Omega$~be a bounded pseudoconvex domain in~$\C^n$ with
  class~$C^\infty$ boundary, and suppose $1\le q\le n$. If a
  compactness estimate~\eqref{compactness estimate} holds for the
  $\dbar$-Neumann problem on $(0,q)$~forms, then the
  $\dbar$-Neumann operator~$N_q$ is a compact (in particular,
  continuous) operator from $W^s_{(0,q)}(\Omega)$ into itself for
  every non-negative~$s$.
\end{theorem}

It suffices to prove the result for integral~$s$, as the
intermediate cases then follow from standard interpolation
theorems \cite{BerghLofstrom1976, PerssonA1964}.  We sketch the
argument for $s=1$, which illustrates the method.  To prove the
compactness of the \(\dbar\)-Neumann operator in
\(W^1_{(0,q)}(\Omega)\), we will establish the (a priori)
estimate \(\|N_qu\|_1^2\le \epsilon \|u\|_1^2+C_\epsilon
\|u\|_0^2\) for arbitrary positive~\(\epsilon\) under the
assumption that \(u\) and \(N_q u\) are both in
\(C^\infty(\overline\Omega)\).

First we show that the compactness estimate~\eqref{compactness
  estimate} lifts to $1$-norms: namely,
\begin{math}
  \|u\|_1^2\le\epsilon(\|\dbar u\|_1^2+\|\dbarstar u\|_1^2)
      +C_\epsilon\|u\|_{-1}^2
\end{math}
for smooth forms~\(u\) in \(\dom\dbarstar\)
(with a new constant~$C_\epsilon$).  In a neighborhood of a
boundary point, we complete $\dbar \rho$ to an orthogonal basis
of $(0,1)$~forms and choose dual vector fields.  To estimate
tangential derivatives of~\(u\), we apply~\eqref{compactness
  estimate} to these derivatives (valid since they preserve the
domain of~$\dbarstar$). We then commute the derivatives with
\(\dbar\) and~\dbarstar, which gives an error term that is of the
same order as the quantity on the left-hand side that we are
trying to estimate, but multiplied by a factor of~\(\epsilon\).
We also need to estimate the normal derivative of~$u$, but since
the boundary is noncharacteristic for the elliptic complex
$\dbar\oplus\dbarstar$, the normal derivative of~$u$ can be
expressed in terms of $\dbar u$, $\dbarstar u$, and tangential
derivatives of~$u$.  Summing over a collection of special
boundary charts that cover the boundary, and using interior
elliptic regularity to estimate the norm on a compact set, we
obtain an inequality of the form \(\|u\|_1^2\le A\epsilon(\|\dbar
u\|_1^2+\|\dbarstar u\|_1^2+\|u\|_1^2) + B(\|\dbar
u\|_0^2+\|\dbarstar u\|_0^2 +\|u\|_0^2)\), where the constants
\(A\) and~\(B\) are independent of~\(\epsilon\). We can use the standard
interpolation inequality \(\|f\|_s\le \epsilon\|f\|_{s+1}
+C_\epsilon \|f\|_{s-1}\) to absorb terms into the left-hand side
when \(\epsilon\)~is sufficiently small.

The lifted compactness
estimate together with the \(L^2\)
boundedness of the \(\dbar\)-Neumann operator implies
\begin{equation}
\|N_qu\|_1^2\le\epsilon(\|\dbar N_qu\|_1^2 +\|\dbarstar
N_qu\|_1^2) +C_\epsilon\|u\|_0^2.\label{one estimate}
\end{equation}
Working as before in special boundary charts, we commute
derivatives and integrate by parts on the right-hand side to make
\(\dbar\dbarstar N_q u +\dbarstar\dbar N_q u=u\) appear (see
\cite[p.~140]{KohnJ1984}, \cite[p.~31]{BoasStraube1990}). Keeping
track of commutator error terms and applying the Cauchy-Schwarz
inequality, we find
\begin{equation}
  \|\dbar N_qu\|_1^2 + \|\dbarstar N_qu\|_1^2 \le A(
  \|N_qu\|_1\|u\|_1 +\|u\|_0^2 + (\|\dbar N_qu\|_1+\|\dbarstar
  N_qu\|_1)\|N_q u\|_1)
\end{equation}
for some constant~\(A\). Consequently \(\|\dbar N_qu\|_1^2 +
\|\dbarstar N_qu\|_1^2 \le B(\|N_qu\|_1^2+\|u\|_1^2)\) for some
constant~\(B\). Combining this with~\eqref{one estimate} gives
the required a priori estimate
\begin{equation}
\|N_qu\|_1^2\le \epsilon\|u\|_1^2 +C_\epsilon\|u\|_0^2.
\end{equation}
Kohn and Nirenberg \cite{KohnNirenberg1965} developed the method
of elliptic regularization (described above after
Theorem~\ref{strictly}) to convert these a priori estimates into
genuine ones.

There is a large class of domains for which the $\dbar$-Neumann
operator is compact \cite{CatlinD1984b, SibonyN1987b}.  In
\cite{CatlinD1984b}, Catlin introduced ``property~(P)'' and
showed that it implies a compactness estimate~\eqref{compactness
  estimate} for the \dbar-Neumann problem. A domain~\(\Omega\)
has property~(P) if for every positive number~\(M\) there exists
a plurisubharmonic function~\(\lambda\) in
\(C^\infty(\overline\Omega)\), bounded between \(0\) and~\(1\),
whose complex Hessian has all its eigenvalues bounded below
by~\(M\) on~\(b\Omega\):
\begin{equation}
  \sum_{j,k=1}^n\frac{\partial^2\lambda}{\partial z_j\partial\bar
    z_k}(z) w_j\bar w_k \ge M|w|^2,\qquad z\in b\Omega, \quad
  w\in\C^n.
\end{equation}
That property~(P) implies a compactness
estimate~\eqref{compactness estimate} follows directly
from~\eqref{basic inequality} and interior elliptic regularity.

It is easy to see that the existence of a strictly
plurisubharmonic defining function implies property~(P), so
strictly pseudoconvex domains satisfy property~(P). So do
pseudoconvex domains of finite type: this was proved by Catlin
\cite{CatlinD1984b} as a consequence of his and D'Angelo's
analysis of finite type boundaries \cite{CatlinD1984a,
  D'AngeloJ1982}.

Property~(P) is, however, much more general than the condition of
finite type. For instance, it is easy to see that a domain that
is strictly pseudoconvex except for one infinitely flat boundary
point must have property~(P).  More generally, property~(P) holds
if the set of weakly pseudoconvex boundary points has Hausdorff
two-dimensional measure equal to zero \cite{BoasH1988,
  SibonyN1987b}.  Sibony \cite{SibonyN1987b} made a systematic
study of the property (under the name of ``B-regularity''). In
particular, he found examples of B-regular domains whose boundary
points of infinite type form a set of positive measure.

It is folklore that an analytic disc in the boundary of a
pseudoconvex domain in~\(\C^2\) obstructs compactness of the
\dbar-Neumann problem: this can be proved by an adaptation of the
argument used in \cite{CatlinD1981, DiederichPflug1981} to show
(in any dimension) that analytic discs in the boundary preclude
hypoellipticity of~\(\dbar\). In higher dimensions, tamely
embedded analytic discs in the boundary obstruct compactness, but
the general situation seems not to be understood; see
\cite{KrantzS1988, LigockaE1985} for a discussion of some
interesting examples.  Salinas found an obstruction to
compactness phrased in terms of the \(C^*\)-algebra generated by
the operators of multiplication by coordinate functions (see the
survey \cite{SalinasN1991} and its references).

In view of the maximum principle, property~(P) excludes analytic
structure from the boundary: in particular, the boundary cannot
contain analytic discs. However, the absence of analytic discs in
the boundary does not guarantee property~(P)
\cite[p.~310]{SibonyN1987b}, although it does in the special
cases of convex domains and complete Reinhardt domains
\cite[Prop.~2.4]{SibonyN1987b}. 

It is not yet understood how much room there is between
property~(P) and compactness.  Having necessary and sufficient
conditions on the boundary of a domain for compactness of the
\dbar-Neumann problem would shed considerable light on the
interactions among complex geometry, pluripotential theory, and
partial differential equations.

\section{The vector field method}
\label{vector field}
In the preceding section, we saw that the \dbar-Neumann problem
is globally regular in domains that support bounded
plurisubharmonic functions with arbitrarily large complex Hessian
at the boundary.  Now we will discuss a method that applies, for
example, to domains admitting defining functions that are
plurisubharmonic on the boundary. The method is based on the
construction of certain vector fields that almost commute
with~\dbar.

We begin with some general remarks about proving a priori
estimates of the form \(\|N_qu\|_s\le C\|u\|_s\) and
\(\|P_qu\|_s\le C\|u\|_s\) in Sobolev spaces for the
\dbar-Neumann operator and the Bergman projection. Firstly, all
the action is near the boundary.  This is clear for the Bergman
projection on functions, because the mean-value property shows
that every Sobolev norm of a holomorphic function on a compact
subset of a domain is dominated by a weak norm on the whole
domain (for instance, the \(L^2\) norm). The corresponding
property holds for the \dbar-Neumann operator due to interior
elliptic regularity.

Secondly, the conjugate holomorphic derivatives
\(\partial/\partial\bar z_j\) are always under control. This is
obvious for the case of the Bergman projection~\(P_0\) on
functions (since holomorphic functions are annihilated by
anti-holomorphic derivatives), and the inequality~\eqref{free}
shows that anti-holomorphic derivatives are tame for the
\dbar-Neumann problem.

Thirdly, differentiation by vector fields whose restrictions to
the boundary lie in the complex tangent space is also innocuous.
Indeed, integrating by parts turns tangential vector fields of
type \((1,0)\) into vector fields of type \((0,1)\), which are
tame, plus lower-order divergence terms
\cite[formula~(3)]{BoasStraube1991}.

Thus, we only need to estimate derivatives in the complex normal
direction near the boundary. Moreover, since the bar derivatives
are free, it will do to estimate either the real part or the
imaginary part of the complex normal derivative. That is, we can
get by with estimating either the real normal derivative, or a
tangential derivative that is transverse to the complex tangent
space.

A simple application of these ideas shows, for example, that the
Bergman projection~\(P_0\) on functions for every bounded
Reinhardt domain~\(\Omega\) in~\(\C^n\) with class~\(C^\infty\)
boundary is continuous from \(W^s(\Omega)\) to itself for every
positive integer~\(s\) \cite{BoasH1984, StraubeE1986}.  Indeed,
the domain is invariant under rotations in each variable, so the
Bergman projection commutes with each angular derivative
\(\partial/\partial\theta_j\). At every boundary point, at least
one of these derivatives is transverse to the complex tangent
space, so \(\|P_0u\|_1\le
C\sum_{j=1}^n\|(\partial/\partial\theta_j)P_0u\|_0 =
C\sum_{j=1}^n \|P_0(\partial u/\partial\theta_j)\|_0 \le
C'\|u\|_1\). Higher derivatives are handled analogously. A similar
technique proves global regularity of the \dbar-Neumann operator
on bounded pseudoconvex  Reinhardt domains
\cite{BoasChenStraube1988, ChenS1989}.

Thus, the nicest situation for proving estimates in Sobolev norms
for the \dbar-Neumann operator is to have a tangential vector
field, transverse to the complex tangent space, that commutes
with the \dbar-Neumann operator, or what is nearly the same
thing, that commutes with \(\dbar\) and~\dbarstar.  (This method
is classical \cite{DerridjM1978, DerridjTartakoff1976,
  KomatsuG1976}.)  Actually, it would be enough for the
commutator with each anti-holomorphic derivative
\(\partial/\partial\bar z_j\) to have vanishing \((1,0)\)
component in the complex normal direction.  However, work of
Derridj \cite[Th{\'e}or{\`e}me~2.6 and the remark following
it]{DerridjM1991} shows that no such field can exist in general.

If we have a real tangential vector field~\(T\), transverse to
the complex tangent space, whose commutator with each
\(\partial/\partial\bar z_j\) has \((1,0)\) component in the
complex normal direction of modulus less than~\(\epsilon\), then
we get an estimate of the form \(\|T^sN_qu\|_0\le A_s(\|u\|_s+
\epsilon\|N_qu\|_s) + C_{s,T}\|u\|_0\).  If the field~\(T\) is
normalized so that its coefficients and its angle with the
complex tangent space are bounded away from zero, then
\(\|T^sN_qu\|_0\) controls \(\|N_qu\|_s\) (independently
of~\(\epsilon\)), so we get global regularity of~\(N_q\) up to a
certain level in the Sobolev scale.  (By making estimates
uniformly on a sequence of interior approximating strongly
pseudoconvex domains, we can convert the a priori estimates to
genuine ones.) Moreover, it suffices if \(T\)~is approximately
tangential in the sense that its normal component is of
order~\(\epsilon\). (This idea comes from work of Barrett
\cite{BarrettD1986}; see the proof of Theorem~\ref{local
  geometry}.)  If we can find a sequence of such normalized
vector fields corresponding to progressively smaller values
of~\(\epsilon\), then the \dbar-Neumann problem is globally
regular at every level in the Sobolev scale.  Because of the
local regularity at points of finite type, the vector fields need
exist only in (progressively smaller) neighborhoods of the
boundary points of infinite type.  In other words, we have the
following result (where the imaginary parts of the~\(X_\epsilon\)
correspond to the vector fields described above)
\cite{BoasStraube1991, BoasStraube1993}.

\begin{theorem}
  \label{vector fields}
  Let $\Omega$ be a bounded pseudoconvex domain in~$\C^n$ with
  class \(C^\infty\) boundary and defining function~\(\rho\).
  Suppose there is a positive constant~$C$ such that for every
  positive~$\epsilon$ there exists a vector field~$X_\epsilon$ of
  type $(1,0)$ whose coefficients are smooth in a
  neighborhood~$U_\epsilon$ in~$\C^n$ of the set of boundary
  points of~$\Omega$ of infinite type and such that
\begin{enumerate}
\item $|\arg X_\epsilon\rho|<\epsilon$ on~$U_\epsilon$, and
  moreover $C^{-1}< |X_\epsilon\rho|<C$ on~$U_\epsilon$,
  and\label{13:1}
\item when $1\le j\le n$,
the form $\partial\rho$ applied to the commutator
$[X_\epsilon,\partial/\partial\bar z_j]$
has modulus less than~$\epsilon$ on~$U_\epsilon$.\label{13:2}
\end{enumerate}
Then the $\dbar$-Neumann operators~$N_q$ (for \(1\le q\le n\))
and the Bergman projections~$P_q$ (for \(0\le q\le n\)) are
continuous on the Sobolev space $W^s_{(0,q)}(\Omega)$ when $s\ge
0$.
\end{theorem}

For a simple example in which the hypothesis of this theorem can
be verified, consider a ball with a cap sliced off by a real
hyperplane, and the edges rounded. The normal direction to the
hyperplane will serve as~\(X_\epsilon\) (the \(U_\epsilon\) being
shrinking neighborhoods of the flat part of the boundary), so the
\dbar-Neumann operator for this domain is continuous at every
level in the Sobolev scale.

Indeed, the hypothesis of Theorem~\ref{vector fields} can be
verified for all convex domains. (The regularity of the
\dbar-Neumann problem for convex domains in dimension two was
obtained independently by Chen \cite{ChenS1991a} using related
ideas.) More generally, the theorem applies to domains admitting
a defining function that is plurisubharmonic on the boundary
\cite{BoasStraube1991}.  We state this as a separate result and
sketch the proof. (Continuity in \(W^{1/2}(\Omega)\) in the
presence of a plurisubharmonic defining function was obtained
earlier by Bonami and Charpentier \cite{BonamiCharpentier1988,
  BonamiCharpentier1990}.)

\begin{theorem}
  \label{plurisubharmonic}
  Let \(\Omega\) be a bounded pseudoconvex domain in~\(\C^n\)
  with class~\(C^\infty\) boundary. Suppose that \(\Omega\) has a
  \(C^\infty\) defining function~\(\rho\) that is
  plurisubharmonic on the boundary:
  \(\sum_{j,k=1}^n(\partial^2\rho/\partial z_j\partial\bar
  z_k)w_j\bar w_k\ge0\) for all \(z\in b\Omega\) and all
  \(w\in\C^n\).  Then for every positive~\(s\) there exists a
  constant~\(C\) such that for all \(u\in W^s_{(0,q)}(\Omega)\)
  we have
\begin{alignat*}{2}
 \|N_qu\|_s &\le C\|u\|_s, & \quad 1&\le q\le n,\\
 \|P_qu\|_s &\le C\|u\|_s, & \quad 0&\le q\le n.
\end{alignat*}
\end{theorem}

Pseudoconvexity says that on the boundary,
\(\sum_{j,k=1}^n(\partial^2\rho/\partial z_j\partial\bar
z_k)w_j\bar w_k\ge0\) for vectors~\(w\) in the complex tangent
space: those vectors for which
\(\sum_{j=1}^n(\partial\rho/\partial z_j)w_j=0\). The hypothesis
of the theorem is that on the boundary, the complex Hessian
of~\(\rho\) is non-negative on all vectors, not just complex
tangent vectors.  (There are examples of pseudoconvex domains,
even with real-analytic boundary, that do not admit such a
defining function even locally \cite{BehrensM1984, BehrensM1985,
  FornaessJ1979}.)  We now sketch how this extra information can
be used to construct the special vector fields needed to invoke
Theorem~\ref{vector fields}.

The key observation is that for each~\(j\), derivatives of
\(\partial\rho/\partial z_j\) of type \((0,1)\) in directions
that lie in the null space of the Levi form must vanish. Indeed,
if \(\partial/\partial\bar z_1\) (say) is in the null space of
the Levi form at a boundary point~\(p\), then
\(\partial^2\rho/\partial z_1\partial\bar z_1(p)=0\), but since
the matrix \(\partial^2\rho/\partial z_j\partial\bar z_k(p)\) is
positive semidefinite, its whole first column must vanish.
(It was earlier observed by Noell \cite{NoellA1991} that the
unit normal to the boundary of a convex domain is constant along
Levi-null curves.)

To construct the required global vector field, it will suffice to
construct a vector field whose commutator with each complex
tangential field of type \((1,0)\) has vanishing component in the
complex normal direction at a specified boundary point~\(p\).
Indeed, these components will be bounded by~\(\epsilon\) in a
neighborhood of~\(p\) by continuity, and we can use a partition
of unity to patch local fields into a global field.  (Terms in
the commutator coming from derivatives of the partition of unity
cause no difficulty because they are complex tangential.)  It is
easy to extend the field from the boundary to the inside of the
domain to prescribe the proper commutator with the complex normal
direction.

Suppose that \(\partial\rho/\partial z_n(p)\ne 0\). We want to
correct the field \((\partial\rho/\partial
z_n)^{-1}(\partial/\partial z_n)\) by subtracting a linear
combination of complex tangential vector fields so as to adjust
the commutators. Since the Levi form may have some zero
eigenvalues at~\(p\), we need a compatibility condition to solve
the resulting linear system. The observation above that type
\((0,1)\) derivatives in Levi-null directions annihilate
\(\partial\rho/\partial z_n\) at~\(p\) is precisely the condition
needed for solvability.  For details of the proof,
see~\cite{BoasStraube1991}.

Kohn \cite{KohnJ1999} has found a new proof and generalization of
Theorem~\ref{plurisubharmonic}. According to a theorem of
Diederich and Forn{\ae}ss \cite{DiederichFornaess1977b} (see also
\cite{RangeR1981}), a smooth bounded pseudoconvex domain admits a
defining function such that some (small) positive power of its
absolute value is plurisuperharmonic inside~\(\Omega\); let
\(\delta\) denote the supremum of such exponents. Kohn showed
that there is a constant~\(A\) such that the \dbar-Neumann
problem is regular in \(W^s(\Omega)\) when \((1-\delta)sA^s<1\).
(This result also contains Proposition~\ref{thm:epsilon}.)

Theorem~\ref{vector fields} applies to other situations besides
the one described in Theorem~\ref{plurisubharmonic}. For
instance, it is possible to construct the vector fields on
pseudoconvex domains that are regular in the sense of Diederich
and Forn{\ae}ss \cite{DiederichFornaess1977c} and Catlin
\cite{CatlinD1984b}. (This gives no new theorem, however, since
the \dbar-Neumann problem is known to be compact on such domains
\cite{CatlinD1984b}; nor does it give a simplified proof of
global regularity in the finite type case, since the construction
of the vector fields still requires Catlin's stratification of
the set of weakly pseudoconvex points \cite{CatlinD1984a}.)

As mentioned in section~\ref{general}, global regularity for the
\dbar-Neumann problem breaks down on the Diederich-Forn{\ae}ss
worm domains.  On those domains, the set of weakly pseudoconvex
boundary points is precisely an annulus, and it is possible to
compute directly that the vector fields specified in
Theorem~\ref{vector fields} cannot exist on this annulus.

For domains of this kind, where the boundary points of infinite
type form a nice submanifold of the boundary, there is a natural
condition that guarantees the existence of the vector fields
needed to apply Theorem~\ref{vector fields}. Following the
notation of \cite{D'AngeloJ1987b, D'AngeloJ1993}, we let
$\eta$~denote a purely imaginary, non-vanishing one-form on the
boundary~$b\Omega$ that annihilates the complex tangent space and
its conjugate.  Let $T$~denote the purely imaginary tangential
vector field on~$b\Omega$ orthogonal to the complex tangent space
and its conjugate and such that $\eta(T)\equiv 1$.  Up to sign,
the Levi form of two complex tangential vector fields $X$ and~$Y$
is $\eta([X,\overline Y\,])$.  The (real) one-form $\alpha$~is
defined to be minus the Lie derivative of~$\eta$ in the direction
of~$T$:
\begin{equation}
\alpha=-\mathcal{L}_T\eta.
\end{equation}
One can show \cite[\S2]{BoasStraube1993} that if \(M\)~is a
submanifold of the boundary whose real tangent space is contained
in the null space of the Levi form, then the restriction of the
form~\(\alpha\) to~\(M\) is closed, and hence represents a
cohomology class in the first De~Rham cohomology \(H^1(M)\). (In
the special case when $M$~is a complex submanifold, this
closedness corresponds to the pluriharmonicity of certain
argument functions, as in \cite{BarrettFornaess1988},
\cite[Prop.~3.1]{BedfordFornaess1978},
\cite[Lemma~1]{BedfordFornaess1981}, and
\cite[p.~290]{DiederichFornaess1977a}.)  This class is
independent of the choice of~\(\eta\).  If this cohomology class
vanishes on such a submanifold~\(M\), and if \(M\)~contains the
points of infinite type, then the vector fields described in
Theorem~\ref{vector fields} do exist.  Thus, we have the
following result~\cite{BoasStraube1993}.

\begin{theorem}
\label{cohomology}
  Let $\Omega$ be a bounded pseudoconvex domain in~$\C^n$ with
  class~\(C^\infty\) boundary.  Suppose there is a smooth real
  submanifold~$M$ (with or without boundary) of~$b\Omega$ that
  contains all the points of infinite type of~$b\Omega$ and whose
  real tangent space at each point is contained in the null space
  of the Levi form at that point (under the usual identification
  of $\R^{2n}$ with $\C^n$).  If the $H^1(M)$ cohomology class
  $[\alpha|_M]$ is zero, then the $\dbar$-Neumann operators~$N_q$
  (for \(1\le q\le n)\) and the Bergman projections~$P_q$ (for
  \(0\le q\le n\)) are continuous on the Sobolev space
  $W^s_{(0,q)}(\Omega)$ when $s\ge0$.
\end{theorem}

On the worm domains, one can compute directly that the class
\([\alpha|_M]\) is not zero.  The appearance of this cohomology
class explains, in particular, why an analytic annulus in the
boundary of the worm domains is bad for Sobolev estimates, while
an annulus in the boundary of other domains may be
innocuous~\cite{BoasStraube1992a}, and an analytic disc is always
benign~\cite{BoasStraube1992a}.  In the special case that \(n=2\)
and \(M\)~is a bordered Riemann surface, Barrett has shown that
there is a pluripolar subset of \(H^1(M)\) such that estimates in
\(W^k(\Omega)\) fail for sufficiently large~\(k\) if
\([\alpha|_M]\) lies outside this subset \cite{BarrettD1997}.
When \(M\)~is a complex submanifold of the boundary,
\([\alpha|_M]\) has a geometric interpretation as a measure of
the winding of the boundary of~\(\Omega\) around~\(M\)
(equivalently, the winding of the vector normal to the boundary).
For details, see \cite{BedfordFornaess1978}. (In the context of
Hartogs domains in~\(\C^2\), see also \cite{BoasStraube1992a}.)

The constructions of the vector fields (needed to apply
Theorem~\ref{vector fields}) in the proofs of Theorems
\ref{plurisubharmonic} and~\ref{cohomology} are more closely
related than appears at first glance.  The vector fields can be
written locally in the form \(e^h L_n +\sum_{j=1}^{n-1} a_j
L_j\), where \(L_1, \dots, L_{n-1}\) form a local basis for the
tangential vector fields of type \((1,0)\), \(L_n\)~is the normal
field of type \((1,0)\), and \(h\) and the \(a_j\) are smooth
functions.  The commutator conditions in Theorem~\ref{vector
  fields} in directions not in the null space of the Levi form
can \emph{always} be satisfied by using the \(a_j\) to correct
the commutators.  Computing the commutators in the remaining
directions leads to the equation
\(dh|_{\mathcal{N}(p)}=\alpha|_{\mathcal{N}(p)}\) at points~\(p\)
of infinite type (where \(\mathcal{N}(p)\) is the null space of
the Levi form at~\(p\)).  The above proof of
Theorem~\ref{plurisubharmonic} amounts to showing that
\(\alpha|_{\mathcal{N}(p)}=0\) when there is a defining function
that is plurisubharmonic on the boundary, whence \(h\equiv0\)
gives a solution.  In Theorem~\ref{cohomology}, the hypothesis of
the vanishing of the cohomology class of~\(\alpha\) on~\(M\)
allows us to solve for~\(h\) (on~\(M\)).

In general, the points of infinite type need not lie in a
``nice'' submanifold of the boundary. It is not known what should
play the role of the cohomology class \([\alpha|_M]\) in the
general situation.  (Note that the analogue of the property that
\(\alpha|_M\) is closed holds in general:
\(d\alpha|_{\mathcal{N}(p)}=0\); see
\cite[\S2]{BoasStraube1993}.)  Furthermore, it is not understood
how to combine the ideas of this section with the pluripotential
theoretic methods discussed in section~5
($B$-regularity/property~(P)).

\section{The Bergman projection on general domains}
In pseudoconvex domains, global regularity of the \dbar-Neumann
problem is essentially equivalent to global regularity of the
Bergman projection \cite{BoasStraube1990}.  In nonpseudoconvex
domains, the \dbar-Neumann operator may not exist, yet the
Bergman projection is still well defined. Since global regularity
of the Bergman projection on functions is intimately connected to
the boundary regularity of biholomorphic and proper holomorphic
mappings \cite{BedfordE1984, BellS1981, BellS1984, BellS1990,
  BellLigocka1980, BellCatlin1982, DiederichFornaess1982,
  ForstnericF1993}, it is interesting to analyze the Bergman
projection directly, without recourse to the \dbar-Neumann
problem. Even very weak regularity properties of the Bergman
projection can be exploited in the study of biholomorphic
mappings \cite{BarrettD1986, LempertL1986}.

In this section, we survey the theory of global regularity of
the Bergman projection on general (that is, not necessarily
pseudoconvex) domains. 

The first regularity results for the Bergman projection that
were obtained without the help of the \dbar-Neumann theory are
in \cite{BellBoas1981}, where it is shown that the Bergman
projection~\(P\) on functions maps the space
\(C^\infty(\overline\Omega)\) of functions smooth up to the
boundary continuously into itself when \(\Omega\)~is a bounded
complete Reinhardt domain with \(C^\infty\)~smooth boundary.

This result was generalized in \cite{BarrettD1982} to domains
with ``transverse symmetries.'' A domain~\(\Omega\) is said to
have transverse symmetries if it admits a Lie group~\(G\) of
holomorphic automorphisms acting transversely in the sense that
the map \(G\times\Omega\to\Omega\) taking \((g,z)\) to \(g(z)\)
extends to a smooth map
\(G\times\overline\Omega\to\overline\Omega\), and for each point
\(z_0\in b\Omega\) the map \(g\mapsto g(z_0)\) of \(G\) to
\(b\Omega\) induces a map on tangent spaces \(T_{\text{Id}}G\to
T_{z_0}^{\R}(b\Omega)\) whose image is not contained in the
complex tangent space to~\(b\Omega\) at~\(z_0\). In other words,
there exists for each boundary point~\(z_0\) a one-parameter
family of automorphisms of~\(\overline\Omega\) whose
infinitesimal generator is transverse to the tangent space
at~\(z_0\). This class of domains includes many Cartan domains as
well as all smooth bounded Reinhardt domains; in both cases,
suitable Lie groups of rotations provide the transverse
symmetries \cite{BarrettD1982}.  For domains with transverse
symmetries, it was observed in \cite{StraubeE1986} that the
Bergman projection not only maps the space
\(C^\infty(\overline\Omega)\) into itself, but actually preserves
the Sobolev spaces.

More generally, one can obtain regularity results in the presence
of a transverse vector field of type \((1,0)\) with holomorphic
coefficients, even if it does not come from a family of
automorphisms.  David Barrett obtained the following result
\cite{BarrettD1986}.

\begin{theorem}
  Let \(\Omega\) be a bounded domain in~\(\C^n\) with
  class~\(C^\infty\) boundary and defining function~\(\rho\).
  Suppose there is a vector field~\(X\) of type \((1,0)\) with
  holomorphic coefficients in \(C^\infty(\overline\Omega)\) that
  is nowhere tangent to the boundary of~\(\Omega\) and such that
  \(|\arg X\rho|<\pi/4k\) for some positive integer~\(k\). Then
  the Bergman projection on functions maps the Sobolev space
  \(W^k(\Omega)\) continuously into itself.
\label{local geometry}
\end{theorem}

In particular, Theorem~\ref{local geometry} implies that there
are no local obstructions to \(W^k\)~regularity of the Bergman
projection. In other words, any sufficiently small piece of
\(C^\infty\) boundary can be a piece of the boundary of a
domain~\(G\) whose Bergman projection is continuous in
\(W^k(G)\): indeed, \(G\)~can be taken to be a small perturbation
of a ball, and then the radial field satisfies the hypothesis of
the theorem.

Theorem~\ref{local geometry} also applies when \(k=1/2\) and the
boundary is only Lipschitz smooth.  For example, the hypothesis
holds for \(k=1/2\) when the domain is strictly star-shaped.
Lempert \cite{LempertL1986} has exploited this weak regularity
property to prove a H\"older regularity theorem for biholomorphic
mappings between star-shaped domains with real-analytic
boundaries.

The first step in the proof of Theorem~\ref{local geometry} is
one we have seen before in section~\ref{vector field}: namely, it
suffices to estimate derivatives of holomorphic functions in a
direction transverse to the boundary. Thus, to bound \(\|Pf\|_k\)
it suffices to bound \(\|X^k Pf\|_0\). However, the inner product
\(\langle X^kPf, X^kPf\rangle\) is bounded above by a constant
times \(|\langle \varphi^k X^k Pf, X^k Pf\rangle|\) when
\(\Re\varphi^k\) is bounded away from zero. By the hypothesis of
the theorem, we can take~\(\varphi\) to be a smooth function that
equals \(\overline X\rho/X\rho\) near the boundary. We then
replace \(\varphi^k X^k\) on the left-hand side of the inner
product by \( (\varphi X-\overline X)^k\), making a lower-order
error (since \(\overline X\)~annihilates holomorphic functions).
The point is that 
\( (\varphi X-\overline X) \) is tangential at the
boundary, so we can integrate by parts without boundary terms,
obtaining \( |\langle Pf, X^{2k}Pf\rangle|\) plus lower-order
terms. Since \(X\)~is a holomorphic field, we can remove the
Bergman projection operator from the left-hand side of the inner
product, integrate by parts, and apply the Cauchy-Schwarz
inequality to get an upper bound of the form
\(C\|f\|_k\|Pf\|_k\). This gives an a priori estimate
\(\|Pf\|_k\le C\|f\|_k\). 
The estimate can be converted into a genuine estimate via an
argument involving the resolvent of the semigroup generated
by the real part of~\(X\). 
For details of the proof, see~\cite{BarrettD1986}.

It is possible to combine such methods with techniques based on
pseudoconvexity. Estimates for the Bergman projection and the
\dbar-Neumann operator on pseudoconvex domains that have
transverse symmetries on the complement of a compact subset of
the boundary consisting of points of finite type were obtained in
\cite{ChenS1987} and~\cite{BoasChenStraube1988}.

A domain in~\(\C^2\) is called a Hartogs domain if, with each of
its points \((z,w)\), it contains the circle \(\{(z,\lambda w):
|\lambda|= 1\}\); it is complete if it also contains the disc
\(\{(z,\lambda w): |\lambda|\le 1\}\).  The (pseudoconvex) worm
domains \cite{DiederichFornaess1977a} and the (nonpseudoconvex)
counterexample domains in \cite{BarrettD1984,
  BarrettFornaess1986} with irregular Bergman projections are
incomplete Hartogs domains in~\(\C^2\).  It is easy to see that
when a Hartogs domain in~\(\C^2\) is complete, the obstruction to
regularity identified in section~\ref{vector field} cannot occur
(see \cite[\S1]{BoasStraube1992a}).  Actually, completeness
guarantees that the Bergman projection is regular whether or not
the domain is pseudoconvex \cite{BoasStraube1989}.  (See
\cite{BoasStraube1992a} for a systematic study of the Bergman
projection on Hartogs domains in~\(\C^2\).)

\begin{theorem}
  Let \(\Omega\) be a bounded complete Hartogs domain
  in~\(\C^2\) with class~\(C^\infty\) boundary. The Bergman
  projection maps the Sobolev space \(W^s(\Omega)\) continuously
  into itself when \(s\ge0\).
\end{theorem}

The proof again uses different arguments on different parts of
the boundary.  An interesting new twist occurs in that the
\dbar-Neumann operator of the envelope of holomorphy of the
domain (which is still a complete Hartogs domain) is exploited.

The Bergman projection is known to preserve the Sobolev spaces
\(W^s(\Omega)\) in all cases in which it is known to preserve the
space \(C^\infty(\overline\Omega)\) of functions smooth up to the
boundary (as is the case for the \dbar-Neumann operator on
pseudoconvex domains).  It is an intriguing question whether or
not this is a general phenomenon.

We now turn to the connection between the regularity theory of
the Bergman projection and the duality theory of holomorphic
function spaces, which originates with Bell \cite{BellS1982b}.
When \(k\)~is an integer, let \(A^k(\Omega)\) denote the subspace
of the Sobolev space \(W^k(\Omega)\) consisting of holomorphic
functions, and let \(A^\infty(\Omega)\) denote the subspace of
\(C^\infty(\overline\Omega)\) consisting of holomorphic
functions. We may view the Fr\'echet space \(A^\infty(\Omega)\)
as the projective limit of the Hilbert spaces \(A^k(\Omega)\),
and we introduce the notation \(A^{-\infty}(\Omega)\) for the
space \(\bigcup_{k=1}^\infty A^{-k}(\Omega)\), provided with the
inductive limit topology.

For discussion of some of the technical properties of these
spaces of holomorphic functions, see \cite{BellBoas1984,
  StraubeE1984}. In particular, \(A^{-\infty}(\Omega)\) is a
Montel space, and subsets of \(A^{-\infty}(\Omega)\) are bounded
if and only if they are contained and bounded in some
\(A^{-k}(\Omega)\). The inductive limit structure on
\(A^{-\infty}(\Omega)\) turns out to be ``nice'' because the
embeddings \(A^{-k}(\Omega)\to A^{-k-1}(\Omega)\) are compact
(as a consequence of Rellich's lemma). Functions in
\(A^{-\infty}(\Omega)\) can be characterized in two equivalent
ways: they have growth near the boundary of~\(\Omega\) that is
at most polynomial in the reciprocal of the distance to the
boundary, and their traces on interior approximating surfaces
\(b\Omega_\epsilon\) converge in the sense of distributions
on~\(b\Omega\).  See \cite{StraubeE1984} for an elementary
discussion of these facts.

The \(L^2\) inner product extends to a more general pairing.
Harmonic functions are a natural setting for this extension. We
use the notations \(h^\infty(\Omega)\) and
\(h^{-\infty}(\Omega)\) for the spaces of harmonic functions
analogous to \(A^{\infty}(\Omega)\) and \(A^{-\infty}(\Omega)\).

\begin{proposition}
\label{duality}
Let \(\Omega\) be a bounded domain in~\(\C^n\) with
class~\(C^\infty\) boundary. For each positive integer~\(k\)
there is a constant~\(C_k\) such that for every
square-integrable harmonic function~\(f\), and every \(g\in
C^\infty(\overline\Omega)\), we have the inequality
\begin{equation}
\label{bell lemma}
\left| \int_\Omega f\bar g\right| \le C_k \|f\|_{-k} \|g\|_k.
\end{equation}
\end{proposition}

The proof of Proposition~\ref{duality} follows from the
observation that for every \(g\in C^\infty(\overline\Omega)\),
there is a function~\(g_1\) vanishing to high order at the
boundary of~\(\Omega\) such that the difference \(g-g_1\) is
orthogonal to the harmonic functions. See 
\cite{BellS1982c}, \cite[Appendix~B]{BoasH1987}, and
\cite{LigockaE1986}
for details; the root idea
originates with Bell \cite{BellS1979} in the context of
holomorphic functions.  Alternatively, Proposition~\ref{duality}
can be derived from elementary facts about the Dirichlet problem
for the Laplace operator \cite{StraubeE1984}.

Because the square-integrable harmonic functions are dense in
\(h^{-\infty}(\Omega)\), it follows from~\eqref{bell lemma} that
the \(L^2\)~pairing extends by continuity to a pairing \(\langle
f,g\rangle\) on \(h^{-\infty}(\Omega)\times
C^\infty(\overline\Omega)\). In particular, this pairing is well
defined and separately continuous on
\(A^{-\infty}_{\mathrm{cl}}(\Omega)\times A^\infty(\Omega)\), where
\(A^{-\infty}_{\mathrm{cl}}(\Omega)\) denotes the closure of
\(A^0(\Omega)\) in \(A^{-\infty}(\Omega)\).

\begin{proposition}
\label{duality equivalence}
Let \(\Omega\) be a bounded domain in~\(\C^n\) with class
\(C^\infty\) boundary. The following statements are equivalent.
\begin{enumerate}
\item The Bergman projection~\(P\) maps the space
\(C^\infty(\overline\Omega)\) continuously into itself.
\item The spaces \(A^{-\infty}_{\mathrm{cl}}(\Omega)\) and
  \(A^\infty(\Omega)\) of holomorphic functions are mutually
  dual via the extended pairing \(\langle\;,\;\rangle\).
\end{enumerate}
\end{proposition}

Proposition~\ref{duality equivalence} is from \cite{BellBoas1984,
  KomatsuG1984}; the case of a strictly pseudoconvex domain is in
\cite{BellS1982b}, and duality of spaces of harmonic functions is
studied in \cite{BellS1982a, LigockaE1986}.  Once
Proposition~\ref{duality} is in hand, Proposition~\ref{duality
  equivalence} is easily proved.  For example, suppose that the
Bergman projection is known to preserve the space
\(C^\infty(\overline\Omega)\), and let \(\tau\) be a continuous
linear functional on the space
\(A^{-\infty}_{\mathrm{cl}}(\Omega)\). Because \(\tau\)~extends
to a continuous linear functional on the inductive limit
\(W^{-\infty}(\Omega)\) of the ordinary Sobolev spaces, it is
represented by pairing with a function~\(g\) in the space
\(W^\infty_0(\Omega)\) of functions vanishing to infinite order
at the boundary. On \(A^0(\Omega)\), and hence on
\(A^{-\infty}_{\mathrm{cl}}(\Omega)\), pairing with~\(g\) is the
same as pairing with~\(Pg\) since, by hypothesis, \(Pg\in
A^\infty(\Omega)\). Therefore \(\tau\)~is indeed represented by
an element of \(A^\infty(\Omega)\).

It is nontrivial that
\(A^{-\infty}_{\mathrm{cl}}(\Omega)=A^{-\infty}(\Omega)\)
when \(\Omega\)~is pseudoconvex. 
Examples show that density properties fail dramatically in the
nonpseudoconvex case \cite{BarrettD1984,
  BarrettFornaess1986}. The arguments in these papers can be
adapted to show that
\(A^{-\infty}_{\mathrm{cl}}(\Omega)\ne A^{-\infty}(\Omega)\) for
these examples.

\begin{theorem}
  Let \(\Omega\) be a bounded pseudoconvex domain in~\(\C^n\)
  with class~\(C^\infty\) boundary. Then the space
  \(A^\infty(\Omega)\) of holomorphic functions is dense both in
  \(A^k(\Omega)\) and in \(A^{-k}(\Omega)\) for each non-negative
  integer~\(k\).
\end{theorem}

The first part is in \cite{CatlinD1980}, the second in
\cite{BellBoas1984}.  In particular, the Bergman projection is
globally regular on a pseudoconvex domain~\(\Omega\) if and only
if the spaces \(A^{-\infty}(\Omega)\) and \(A^\infty(\Omega)\)
are mutually dual via the pairing \(\langle\;,\;\rangle\).

Here is a typical application of Proposition~\ref{duality
  equivalence} to the theory of the Bergman kernel function
\(K(w,z)\).

\begin{corollary}
\label{determinacy}
Let \(\Omega\) be a bounded domain in~\(\C^n\) with
class~\(C^\infty\) boundary. Suppose that the Bergman projection
maps the space \(C^\infty(\overline\Omega)\) into itself. If
\(S\) is a set of determinacy for holomorphic functions
on~\(\Omega\), then \(\{K(\cdot\,,z): z\in S\}\) has dense linear
span in~\(A^\infty(\Omega)\).
\end{corollary}

Indeed, global regularity of the Bergman projection~\(P\)
implies that \(K(\cdot\,,z)\in A^\infty(\Omega)\) for each \(z\)
in~\(\Omega\), since \(K(\cdot\,,z)\) is the projection of a
smooth, radially symmetric bump function (this idea originates
in \cite{KerzmanN1972}). Now if a linear functional~\(\tau\) on
\(A^\infty(\Omega)\) vanishes on each \(K(\cdot\,,z)\) for \(z\in
S\), then \(\tau(z)=0\) on~\(S\), whence \(\tau\equiv0\). (Note
that since \(\tau\in A^{-\infty}_{\mathrm{cl}}(\Omega)\), the
Bergman kernel does reproduce~\(\tau\), because evaluation at an
interior point is continuous in the topology of
\(A^{-\infty}(\Omega)\).)  

Corollary~\ref{determinacy} is due to Bell \cite{BellS1979,
  BellS1982b}. It is the key to certain non-vanishing properties
of the Bergman kernel function that are essential in the approach
to boundary regularity of holomorphic mappings developed by
Bell, Ligocka, and Webster \cite{BellLigocka1980, BellS1981,
  BellS1984, LigockaE1980, LigockaE1981, WebsterS1979}.

\bibliographystyle{amsplain}

%\bibliography{survey}

%% bibliography generated by BiBTeX
\providecommand{\bysame}{\leavevmode\hbox to3em{\hrulefill}\thinspace}

\end{document}